\documentclass[a4paper, 11pt]{article}
\usepackage[T1]{fontenc}
\usepackage{textcomp}
\usepackage{amsmath}
\usepackage{amsfonts}
\usepackage{amssymb}
\usepackage{amsthm}
\usepackage{tikz}
\usetikzlibrary{cd}
\usepackage{cancel}
\usepackage{multicol}
\usepackage{geometry} 
\usepackage{parskip}
\usepackage{multicol}
\usepackage{hyperref}
\usepackage{xcolor}
\usepackage[polish,english]{babel}
\usepackage{csquotes}

\newtheorem{thm}{Theorem}[subsection]
\newtheorem{lm}[thm]{Lemma}
\newtheorem{fact}[thm]{Fact}

\theoremstyle{definition}
\newtheorem{df}[thm]{Definition}
\renewcommand{\thethm}{%
  \ifnum\value{subsection}=0 
    \thesection
  \else
    \thesubsection
  \fi%
  .\arabic{thm}}

\newcommand{\N}{\mathbb{N}}

\newcommand{\R}{\mathbb{R}}

\newcommand{\s}{\mathfrak{s}}
\newcommand{\p}{\mathfrak{p}}
\newcommand{\A}{\mathcal{A}}
\newcommand{\M}{\mathcal{M}}
\newcommand{\omom}{[\omega]^\omega}

\DeclareMathOperator{\cl}{cl}
\let\int\relax
\DeclareMathOperator{\int}{int}

\DeclareMathOperator{\cov}{cov}
\DeclareMathOperator{\Clop}{Clop}
\newcommand{\cone}[1]{\ensuremath{\langle #1 \rangle}}

\let\epsilon\varepsilon
\let\phi\varphi
\let\emptyset\varnothing
\let\le\leqslant

\let\ge\geqslant

\let\bd\partial

\let \b \textbf
\usepackage[backend=biber, sorting=ynt]{biblatex}
\title{Convergent sequences in various topological spaces}
\author{Dawid Migacz}
\date{}
\begin{document}
\newpage
\thispagestyle{empty}
\newgeometry{margin=1in}
\begin{multicols}{2}
\begin{center}
\textbf{\Large University of Wrocław \\}
\vspace{0.15cm}
\textbf{\Large Faculty of Mathematics and~Computer Science \\}
\vspace{0.3cm}
\textbf{\Large Mathematical Institute} \\ \ \\
\textit{\large speciality: theoretical mathematics}
\end{center}
\columnbreak
\begin{center}
\textbf{\Large Uniwersytet Wrocławski \\}
\vspace{0.15cm}
\textbf{\Large Wydział Matematyki i~Informatyki \\}
\vspace{0.3cm}
\textbf{\Large Instytut Matematyczny} \\ \ \\
\textit{\large specjalność: matematyka teoretyczna}
\end{center}
\end{multicols}
\vspace{4cm}
\begin{center}
\textbf{\textit{\Large Dawid Migacz}}\\
\vspace{0.5cm}
\textbf{\LARGE Convergent sequences in various topological spaces}

\textbf{ Ciągi zbieżne w przestrzeniach topologicznych}
\end{center}
\vspace{3cm}
\begin{multicols}{2} \large
\begin{center}
Bachelor's thesis \\ written under the supervision of \\ dr hab. Piotr Borodulin-Nadzieja \end{center}
\columnbreak
\begin{center}
Praca licencjacka \\ napisana pod kierunkiem \\ dr. hab. Piotra Borodulina-Nadziei
\end{center}
\end{multicols}

\vfill
\begin{center}
{\Large Wrocław, 2021}\\
\end{center}
\restoregeometry

\newpage

\setcounter{page}{2}
\begin{abstract}
    The following paper is inspired by Efimov's problem -- an undecided problem of whether there exists an infinite compact topological space that does not contain neither non-trivial convergent sequences nor a copy of $\beta\omega$. After introducing the basic topological concepts,  we present several classes of topological spaces in which such sequences can certainly be found, namely ordered, scattered, metrisable spaces and Valdivia compacta. We show that some cardinal coefficients set limits on the smallest cardinality of the base and the smallest cardinality of a neighbourhood base, under which the existence of convergent sequences can be ensured. In the final part of the paper we define the space $\beta\omega$ and show its selected properties. In particular, we prove that there are indeed no non-trivial convergent sequences in $\beta\omega$. 
\end{abstract}
\begin{otherlanguage}{polish}
\begin{abstract}
Poniższa praca inspirowana jest problemem Jefimowa -- nierozstrzygniętym zagadnieniem, czy każda nieskończona zwarta przestrzeń topologiczna niezawierająca nietrywialnych ciągów zbieżnych zawiera kopię przestrzeni $\beta\omega$. Po wprowadzeniu podstawowych pojęć topologicznych prezentujemy kilka klas przestrzeni topologicznych, w których nietrywialne ciągi zbieżne na pewno się znajdują, mianowicie przestrzenie uporządkowane, rozproszone, metryzowalne oraz kompakty Valdivii.  Pokazujemy, że niektóre współczynniki kardynalne zadają ograniczenia na najmniejszą moc bazy i najmniejszą moc bazy otoczeń, stanowiące warunek wystarczający istnienia ciągów zbieżnych. W końcowej części pracy definiujemy przestrzeń $\beta\omega$ oraz pokazujemy wybrane jej własności. W szczególności dowodzimy, że istotnie w $\beta\omega$ nie ma nietrywialnych ciągów zbieżnych. 
\end{abstract}
\end{otherlanguage}

\tableofcontents
\setcounter{section}{-1}
\section{Introduction}
The study of convergent sequences in compact spaces is an active branch of contemporary topology. However, there is no survey article in the literature about the reasons for their existence. 

We present several topological properties implying that the space having this property contains a non-trivial convergent sequence (i.e. a sequence which is not eventually constant). In particular those classes of compact infinite spaces have non-trivial convergent sequences:
\begin{itemize}
    \item scattered spaces (Theorem \ref{scatteredconv})
    \item ordered spaces (Theorem \ref{orderedconv})
    \item metrisable spaces (Theorem \ref{metrisableconv})
    \item spaces of small character (Theorems \ref{firstcontconv} and \ref{pconv})
    \item spaces of small weight (Theorem \ref{covMconv})
    \item Valdivia compacta (Theorem \ref{valdiviaconv})
\end{itemize}
Notice that this list is not irreducible (e.g. compact metrisable spaces are ordered, they are also Valdivia). Nevertheless, we enclose separate proofs for all those classes, in each case trying to grasp the precise argument.

Whereas the statements of these theorems are commonly known, the proofs are notoriously difficult to find. In this paper we intend to fill that gap.

At the end we present the flagship compact infinite space without non-trivial convergent sequences, namely $\beta\omega$. Note that the problem if every infinite compact space contains either a non-trivial convergent sequence or a copy of $\beta\omega$ is still unsolved (at least in $\mathsf{ZFC}$; it is called the Efimov Problem). 

The proofs of Theorems \ref{scatteredconv}, \ref{covMconv}, \ref{sigmasubsetslimitpoints}, \ref{valdiviaconv} are adaptations of proofs from the literature. I proved the remaining claims myself with the invaluable guidance of my supervisor.

\section{Basic topological concepts}
\begin{df}
A \textbf{topological space} is a set (of points) $X$ together with a collection of some of its subsets. We call these distinguished sets \textbf{open sets}. We call their complements \textbf{closed sets}. We call a set that is simultaneously closed and open \textbf{clopen}. The family of all open sets will be called a \textbf{topology}. In a topological space the following axioms must be satisfied: 
\begin{enumerate}
    \item $X$ and $\emptyset$ are open;
    \item the union of any collection of open sets is open;
    \item the intersection of finitely many open sets is open.
\end{enumerate}
\end{df}
When speaking in context of a space, a set is a subset of said space.
\begin{df}
We say that $x$ is \textbf{isolated} if $\{x\}$ is open.
\end{df}
\begin{df}
The \textbf{interior} of a set $A$ is the $\subseteq$-biggest open set contained in it. Its \textbf{closure} is the $\subseteq$-smallest closed set containing it. We denote them by $\int A$ and $\cl A$ respectively. The \textbf{boundary} of $A$ is defined as $\cl A \setminus \int A$ and denoted by $\bd A$.
\end{df}

\begin{df}
A topological space is \textbf{Hausdorff} if every two points can be separated by open sets, i.e. for all $x,y \in X$ there are disjoint open sets $U, V$, such that $x\in U$ and $y \in V$. 
\end{df}
\begin{df}
We call a topological space \textbf{compact} if every open cover has a finite subcover. \end{df}
Specifically, in a compact space, given an arbitrary (possibly infinite) family of open sets $\{U_\alpha: \alpha < \kappa\}$, such that $\bigcup U_\alpha = X $, we are always able to choose a finite collection \{$U_{a_1},\ldots, U_{a_n}$\}, whose union will still be the whole space: $\bigcup_{k=0}^n U_{a_k}=X$. 

\begin{df}
We define the \textbf{subspace topology} of $Y \subseteq X$ in the following way: a set $V\subseteq Y$ is open in $Y$ only if there exists some $U$ open in $X$ such that $V = Y \cap U$.
\end{df}
\begin{lm} \label{closediscompact}

A closed subspace $Y$ of a compact space $X$ is compact. 
\end{lm}
\begin{proof}
Take a cover $\{V_\alpha\}$ of $Y$. Let $U_\alpha$ be any open set in $X$ such that $V_\alpha = Y \cap U_\alpha$. The family $U_\alpha$, together with $X\setminus Y$, constitutes an open cover of $X$. From compactness of $X$ we can take a subcover: $\{U_{a_1},\ldots,U_{a_n}, (X \setminus Y)\}$. But then, $\{U_{a_1} \cap Y, \ldots, U_{a_n}\cap Y, (X \setminus Y)\cap Y\} = \{V_{a_1},\ldots,V_{a_n}, \emptyset\}$ is, after discarding the empty set, a finite subcover of $Y$.
\end{proof}

\begin{lm}
In a compact Hausdorff space a point $y$ can be separated from a closed set $K$ by open disjoint sets $U \ni x, V \supseteq K$. 
 \label{separClosed}
\end{lm}
\begin{proof}
Note that $K$ must be compact, from the previous lemma. For each $x \in K$ separate it from $y$ by $V_x \ni y$ and $U_x \ni x$. The sets $U_x$ constitute an open cover of $K$, thus there exists a finite open subcover $\{U_{x_n}: n < N\}$. Their union is an open superset of $K$, disjoint with the open set $\bigcap_{n<N} V_{x_n} \ni y$. 
\end{proof}
\begin{df}
We say that a space is \textbf{discrete} (or has a discrete topology) if every set is open.
\end{df}
Clearly, no infinite discrete space can be compact, because the family of all singletons constitutes a cover without any subcovers, let alone finite ones. 

The discrete topology is the natural topology on the set of natural numbers $\N$. From now on, we will denote that set by $\omega$.

\begin{df}
A \textbf{base} of a topological space is a family of sets such that every open set can be written as a union of those base sets. 
\end{df}
\begin{df}
The \textbf{weight} of a topological space is the minimal cardinality of a base.
\end{df}
\begin{df}
We say that a sequence $\cone{a_n}$ is \textbf{convergent} to $a$ if for every open neighbourhood $U$ of $a$ all but finitely many elements of $\cone{a_n}$ are in $U$. We then write $a_n \rightarrow a$. A sequence is \textbf{trivial} if it is eventually constant. 
\end{df}

\section{Spaces containing non-trivial convergent sequences}
Now we proceed to describe several spaces which do have a convergent sequence. We presume all spaces to be compact, infinite and Hausdorff.
\subsection{Scattered spaces}
\begin{df}
We call a space \textbf{scattered} if every closed subspace has an isolated point. 
\end{df}
Possibly the most natural compact infinite example of such a space would be $\omega \cup \{\infty\}$ thought of as the one-point compactification of the discrete space $\omega$, in which all subsets of $\omega$ are open, but a set containing $\infty$ is open iff it is cofinite. 
\begin{thm}
Every scattered space has a convergent sequence.
\label{scatteredconv}
\end{thm}
\begin{proof}
We present a slightly modified version of the proof from \cite{backe}. We are going to construct such a sequence. Let $X$ be an infinite scattered space. Our approach will be to find a subspace similar to the previously mentioned $\omega \cup \{\infty\}$, i.e. a countable closed space, in which almost all singletons are open. We will temporarily forgo the requirement of closedness, only taking closure later. So we are looking for an analogue to $\omega$ -- a space that is countable and discrete. Note that it implies openness. It can be found in two ways.

The first is to use Lemma \ref{discreteinhausdorff}. Having an infinite discrete subspace provided by that lemma, we take one of its countable subspaces and call it $A$. 

The second way is to construct it explicitly. Because $X$ is infinite, closed, and scattered, there is an isolated $p_0 \in X$. Note that $X \setminus \{p_0\} = X \setminus \bigcup_{n<1} \{p_n\}$, as a complement of an open set, is closed, so there exists an isolated $p_1 \in X \setminus \{p_0\}$. It means that $\{p_0,p_1\}$ or $\{p_0\}$ is open. But in compact Hausdorff spaces finite sets are always closed, so we know that $\{p_1\}$, being the intersection of $X \setminus \{p_0\}$ with one of said sets, is indeed open, so $p_1$ is isolated. We can repeat this construction, obtaining an infinite sequence of isolated points, and therefore the space we wanted, $A:=\bigcup_{n<\omega} \{p_n\}$.

Observe that the boundary $\bd A$ of $A$ is non-empty. Otherwise $A$ itself would have to be closed, and as a closed subspace of a compact space, compact. But then it would be impossible to find a finite subcover of a cover consisting of singletons. What's more, the boundary is equal to $ (X \setminus {A}) \cap \cl{A}$, so it is closed itself. Because $A$ is open, no element of $\bd A$ is in $A$. Notice that there are no open sets contained in the boundary, because otherwise the closure would be smaller.

This means that there is an isolated $b \in \bd A$. We can separate it, by Lemma \ref{separClosed}, from the remaining $\bd A \setminus \{b\}$ by open U and R, respectively. 

Let $\cone{b_n : n < \omega}$ be a sequence of all elements of the countable set $B:=A \cap U$. We will show that it is convergent to $b$. First, we claim that $B \cup \{b\}$ is closed, while $B$ is not. 

\begin{enumerate}
    \item Note that for every open neighbourhood $N$ of $b$ the sets $N \cap B$ and $N \setminus (U \cap B)$ are non-empty. The latter contains $b$. The former cannot be empty, because if $N$ was disjoint with the open set $B=U\cap A$, the open intersection $N\cap U$ would have to be contained in $\bd A \setminus A$, which cannot happen, as previously mentioned. Consequently $b \not\in B$ while $b \in \bd B$, so $B$ is not closed.
    \item Of course then $B \cup \{b\} \subseteq \cl B$. We will prove that the equality holds. Suppose $x \in \cl B \setminus B$. Then for every open neighbourhood $N \ni x$ the sets $A\cap (N \cap U) = N \cap B$ and $N \setminus B$ are non-empty. Considering an open neighbourhood $N \cap U$ allows us to conclude that $(N\cap U) \setminus B = (N\cap U) \setminus (A \cap U) = (N\cap U) \setminus A$ is non-empty. Therefore $N\cap A$ and $N \setminus A$ is non-empty, so $x\in\bd A$.
    
    Note that $B$ is a subset of a closed set $X\setminus R$, therefore $\cl B$ is too. Hence $x$ is not in $R$, but the only element of $\bd A$ not in $R$ is $b$, so $x=b$. This asserts that the equality indeed holds.
\end{enumerate}
We see that $B\cup\{b\}$ is closed and thus compact. We will now show that the defined earlier sequence converges to $b$. Take an open neighbourhood $N$ of $b$ in $X$. The restriction $N\cap \cl B$ to $\cl B$ must be open in $\cl B$. If infinitely many points of $B$ (of which all are isolated) were outside $N\cap \cl B$, they would constitute an open cover without a finite subcover, which is a contradiction with compactness. Therefore for every open neighbourhood of $b$ all but finitely many elements of the defined sequence are in that neighbourhood.
\end{proof}

\subsection{Ordered spaces}
\begin{df}
An \textbf{ordered space} is a topological space with a linear ordering, such that the family of all open intervals $(a,b):=\{x\in X:a<x<b\}$ and rays $(-\infty,b):=\{x\in X:x<b\}, (a,\infty):=\{x\in X:a<x\},$  constitutes a base of its topology. 
\end{df}
Let $\cone{a_n: n < \omega}$ be a strictly increasing sequence. Of course it is non-trivial. We will prove that it is convergent. First we need to find its supremum.
\begin{lm}
There exists an upper bound $u$ such that $u>a_n$ for all $n$. 
\end{lm}
\begin{proof}
  Otherwise $\{(-\infty, a_0)\} \cup \{(a_n, a_{n+2}):n<\omega \}$ would constitute an open cover without a finite subcover, which is impossible in compact spaces.
\end{proof}
\begin{lm}
There exists a supremum of $\cone{a_n}$.
\end{lm}
\begin{proof}
For each upper bound $u_\alpha$ take the ray $U_\alpha=(u_\alpha,\infty)$. Suppose none of the upper bounds is the smallest. Then every upper bound (i.e. every point greater than all elements of $\cone{a_n}$) is in some ray $U_\alpha$. Note that the elements that are lesser or equal to at least one $a_n$ are in $\{(-\infty, a_1)\} \cup \{(a_n, a_{n+2}):n<\omega \}$.  But then  $\{(-\infty, a_0)\} \cup \{(a_n, a_{n+2}):n<\omega \} \cup \{U_\alpha, \alpha<\kappa\}$ is an infinite open cover without a finite subcover, which is a contradiction.
\end{proof}
\begin{thm}
There is a non-trivial convergent sequence in an infinite ordered space. 
\label{orderedconv}
\end{thm}
\begin{proof}
Let $a$ be the supremum of $\cone{a_n}$. Take an open interval $(l,r)$ such that $a\in(l,r)$. Because a is the supremum, $l$ cannot be greater than all of the elements of the sequence. It must then be lesser or equal to some $a_n$. But then, all but finitely many elements of the sequence are in $(l,r)$. Therefore $\cone{a_n}$ converges to $a$.
\end{proof}
\subsection{Metrisable spaces}
\begin{df}
A \textbf{metric} is a function $d: X\times X \to [0,\infty)$ such that for all $x,y,z\in X$ the following conditions are met:
\begin{enumerate}
    \item $d(x,y)=0 \iff x=y$,
    \item $d(x,y)=d(y,x),$
    \item $d(x,y)+d(y,z)\ge d(x,z).$
\end{enumerate}
We call the pair $\cone{X,d}$ a \textbf{metric space}. There is a topology connected with a metric space, namely the topology in which all open balls $B_r(x_0):=\{x\in X: d(x,x_0)<r\}$ form a base. We assume $B_\infty(x_0)=X$.
\end{df}

\begin{df}
A topological space is \textbf{metrisable} if there exists a metric on it that generates the topology of that space.   
\end{df}
\begin{fact}
Metrisable spaces are Hausdorff. 
\end{fact}
\begin{proof}
Let $X$ be a metrisable space with a metric $d$. We want to separate $x,y \in X$. They can be separated by the obviously disjoint $B_{\frac13d(x,y)}(x)$ and $B_{\frac13d(x,y)}(y)$.
\end{proof}
\begin{thm}
Every metrisable space contains a non-trivial convergent sequence.
\label{metrisableconv}
\end{thm}
\begin{proof}
Let $X$ be a metrisable space with a metric $d$. If all points were isolated, the space could not be compact. Therefore we can take a non-isolated point $x\in X$. Of course $\bigcap_{n=0}^\infty B_\frac1n(x)=\{x\}$. But finite intersections $\bigcap_{n=0}^N B_\frac1n(x)= B_\frac1N(x)$ are open, thus are not $\{x\}$. Define $A_n = B_\frac1n(x)\setminus B_\frac1{n+1}(x)$. Note $\bigcup_{n\in\N} A_n = X \setminus\{x\}$. The sequence $A_n$ cannot stabilise, because it would imply that $x$ is isolated. Taking a sequence $C_n$ of these sets $A_k$ that are not equal to their predecessors we can choose $c_n \in C_n$ which are convergent to $x$. Indeed, $d(c_n,x)<\frac1n$.
\end{proof}
\begin{fact}
A metrisable space is of countable weight.
\end{fact}
\begin{proof}
For every radius $\frac1n, n\in\N$ we see by compactness that finitely many open balls $B_\frac1n (x^n_1), \ldots,  B_\frac1n (x^n_k)$ of that radius cover the whole space. We claim that the family of all such balls constitutes a (obviously countable) base. Take a point $x$ in an arbitrary open set $U$
. Because $U$ is open, there is an $r>0$ such that $B_r(x)\subseteq U$. Take $n$ such that $\frac1n < \frac{r}{2}$. We already know, that for some $l$ we have $x\in B_\frac1n (x^n_l)$. It is obvious that the distance between any two points in such a ball is lesser than $\frac2n$ and therefore lesser than $r$. Hence $x\in B_\frac1n (x^n_l) \subseteq U$; the conclusion follows. 
\end{proof}

We see that metrisable spaces have the smallest possible weight.

\subsection{First-countable spaces}
We can generalise our result from metrisable spaces in a different way, first noting that in such a space every neighbourhood of a point contains a ball centred in that point. 
\begin{df}
A \textbf{neighbourhood base} of $x\in X$ is a collection $\cone{N_\alpha}$ of (some) neighbourhoods of $x$ such that for every neighbourhood $M \ni x$ there is an $\alpha$ such that $N_\alpha \subseteq M$.
\end{df}
\begin{df}
A space is \textbf{first-countable} if every point has a countable neighbourhood base.
\end{df}
\begin{df}
The \b{character} of a point of a topological space  is the minimal cardinality of a neighbourhood base of that point. 
\end{df}
\begin{df}
The \b{character} of a topological space is the supremum of the characters of its points. 
\end{df}
Therefore first-countability of a space is equivalent to that space having countable character. 
\begin{thm} \label{firstcontconv}
All first-countable spaces contain a non-trivial convergent sequence.
\end{thm}
\begin{proof}
Take a non-isolated $x\in X$ and its countable neighbourhood base $\cone{N_n}$. Consider a sequence of open sets $M_n:=N_1\cap\ldots\cap N_n$. First, note that none of these sets is a singleton of $x$, because we demanded $x$ to be non-isolated. 

We will now show that the sequence does not stabilise. Take $y\in M_n \setminus \{x\}$ and separate $x$ and $y$ by open sets $U$ and $V$. But then there is some $N_k \subseteq U$, so $M_k \subset M_n$. Therefore, the sequence $\cone{M_n}$ contains infinitely many different sets. Renumerate them in such a way that that sequence is strictly decreasing. Notice that $\cone{M_n}$ is a neighbourhood base -- for every neighbourhood $U$ of $x$ some $N_n$ is contained in $U$, but then $M_n$ also is. 

Let $x_n \in X_{n} \setminus X_{n+1}$ for every $n$. We will now show that $x_n \rightarrow x$. 
Take an open set $U\ni x$. It must contain a set $M_n$ for some $n$. But then all $x_k$ for $k \ge n$ are in $M_n$ and therefore in $U$. Hence only finitely many elements of the sequence $\cone{x_n}$ lie outside any given open set containing $x$.
\end{proof}
Note that our assumption need not be so strong -- a countable neighbourhood base existing just for one non-isolated point would suffice. Furthermore, we can permit the neighbourhood base to be somewhat larger.

\subsection{Spaces of weight smaller than \texorpdfstring{$\s$}{s} }
In this subsection all sets will be subsets of $\omega$, unless specified otherwise. The family of all infinite subsets of $\omega$ will be denoted $\omom$.
\begin{df}
Given two infinite sets $A,B$ we say that $A$ \textbf{splits} $B$ if both $B \cap A$ and $B \setminus A$ are infinite. 
\end{df}
\begin{df}
Consider a family $\A=\{A_\alpha: \alpha < \kappa \} $ of infinite subsets of $\N$. We say that $\A$ \textbf{is splitting} if for every $B \in \omom$ there is a set $A \in \A$ which splits $B$.
\end{df}
\begin{df}
The \textbf{splitting number} $\s$ is the cardinality of the smallest splitting family. 
\end{df}
\begin{fact}
The splitting number is uncountable.
\end{fact}
\begin{proof}
Take a family $\A=\{A_\alpha: \alpha < \omega \} $. We will recursively will construct a sequence $B_n$ of infinite sets, whose limit (in some sense) will not be split by $\A$. We assume $B_{-1}=\omega$.
The set $B_{-1}\setminus A_0$ must be infinite or have an infinite complement. Let $B_0$ be the infinite one of them. We repeat this operation, setting $B_n$ to be either $B_{n-1} \cap A_n$ or $B_{n-1} \setminus A_n$, depending on which one is infinite. Both cannot be finite, as their union is an infinite set $B_{n-1}$.
Note that  is for every set $A_{n}$ all the following sets $B_{n+k}$ are either disjoint with $A_n$, or contained in it.

Denote the $k$-th element of $B_n$ by $B_n(k)$. Let $C_n=B_n(n)$ and $C=\{C_n:n<\omega\}$. Observe that $\cone{C_n}$ is a strictly increasing sequence. Fix $A_p$. We know that  $C_n\in A_p$ for all $n\ge p$ or  $C_n\not\in A_p$ for all $n\ge p$. Therefore only the sets $C_n$ for $n<p$ can be different in this aspect. So we have proved that either $A_p\setminus C$ is finite or $A_p \cap C$ is finite.

Therefore we have demonstrated that there exists a set not split by any element of $\A$.
\end{proof}
\begin{thm}
All spaces of weight smaller than $\s$ contain a non-trivial convergent sequence.
\label{sconv}
\end{thm}
\begin{proof}

Take a countable sequence $\cone{y_n}$. Note that because all splitting families have cardinality of at least $\s$,  there exists an infinite subsequence $\cone{y_{a_n}}=:\cone{z_n}$ such that no element of the base splits it. Therefore the base can be divided into two classes -- the class $\cone{F_\alpha}$ of the base sets containing only finitely many elements of $\cone{z_n}$ and the class $\cone{A_\alpha}$ of the base sets containing all but finitely many elements of $\cone{z_n}$. Suppose $\cone{F_\alpha}$ is an open cover of $X$. Compactness means that there exists a finite subcover. But then, only finitely many elements of $\cone{z_n}$ can be covered. This contradiction means that there exists a point $z\in X$ not contained in any $F_\alpha$. Therefore, its every neighbourhood has an $A_\alpha$ included, hence every open neighbourhood of $z$ contains all but finitely many elements of $\cone{z_n}$.

\end{proof}

\subsection{Spaces of character smaller than \texorpdfstring{$\p$}{p}}
\begin{df}
We say that $P \subseteq \omega$ is a \textbf{pseudo-intersection} of a family $N$ of subsets of $\omega$ if for every $N_\alpha \in N$ the set $P \setminus F$ is finite.
\end{df}
\begin{df}
We say that a family $N$ of sets has the \textbf{finite intersection property} or, equivalently, is a \textbf{filter base} if every finite subfamily has a non-empty intersection.
\end{df}
\begin{df}
The \textbf{pseudo-intersection number} $\p$ is the cardinality of the smallest filter base not containing an infinite pseudo-intersection. 
\end{df}
\begin{thm}
All spaces containing a smaller than $\p$ neighbourhood base of a non-isolated point $x \in X$  non-trivial convergent sequence.
\label{pconv}
\end{thm}
\begin{proof}

Recall that an infinite space $K$ must have a countable discrete subspace $A$. It cannot be closed in $K$, because a closed subspace of a compact space is compact and there are no infinite compact spaces. Therefore $\bd A \neq A$. Take $x \in \bd A$.
Let $\{N_\alpha:\alpha<\kappa<\p\}$ be a neighbourhood base of $x$. Note that a neighbourhood base of $x$ is a filter base -- finite intersections of base elements are open sets containing $x$. 

Let $M_\alpha := N_\alpha \cap A$. Notice that every neighbourhood of $x$ contains an element of $A$ -- otherwise $x$ would not be in $\cl A$. Therefore the family $\{M_\alpha:\alpha < \kappa < \p\}$ is a filter base on $A$. Hence, there exists a pseudo-intersection $P \subseteq A$ for that family.

Now take an arbitrary neighbourhood $U\ni x$. For some $\alpha$ we have $M_\alpha \subseteq N_\alpha \subseteq U$. But then $P\setminus M_\alpha$ is finite. We have proved that only finitely many elements of $P$ are outside of an arbitrary neighbourhood of $x$, so any enumeration of $P$ works.

\end{proof}

\subsection{Spaces of weight smaller than \texorpdfstring{$\cov(\mathcal{M})$}{cov(M)}}
Now we are going to provide an unnecessary, but nonetheless aesthetically pleasing definition of an object dual to filters -- a family of sets, in some sense, small.
\begin{df}
An \b{ideal} over $X$ is a family $I$ of subsets of $X$ such that:
\begin{enumerate}
    \item $\emptyset \in I$, $X \not\in I$;
    \item if $A, B \in I$, then $A \cup B \in I$;
    \item if $A \in I$ and $B\subseteq A$, then $B \in I$.
    \end{enumerate}
\end{df}
\begin{df}
We call a set $B$ \b{nowhere dense} if $\int \cl B = \emptyset$. 
\end{df}
\begin{df}
We call a set \b{meagre} if it is a countable union of nowhere dense sets.
\end{df}
\begin{df}
The \b{Cantor space} is the metric space $2^\omega$ (the set of all infinite sequences of zeroes and ones), where the distance between two sequences is equal to $\frac{1}{2^n}$, where $n$ is the number of the first place where the sequences differ. 
\end{df}

Note that the balls in that space consist of sequences having a common beginning. Recall that in metric spaces all the balls constitute a canonical base. 
\begin{fact}
The Cantor space has a clopen countable base.
\end{fact}
\begin{proof}
We will show on an example that all balls are closed. Take the ball containing sequences beginning with $\cone{0,1}$, denoted by $B_{\cone{0,1}}$. Its complement is the union of $B_{\cone{0,0}}$, $B_{\cone{1,0}}$ and $B_{\cone{1,1}}$, therefore it is open. 

In general, it is clear that the complement of a ball is the union of finitely many balls. 

The countability is trivial -- the canonical base is a countable union of finite sets -- for every $n<\omega$ we take balls corresponding to the beginnings of length $n$.
\end{proof}

\begin{lm}
A compact space containing a countable clopen base contains countably many clopen sets. 
\end{lm}
\begin{proof}
Let $C$ be an arbitrary clopen set. It has to be equal to an arbitrary union of balls. But $C$ is also closed, so it is compact, by \ref{closediscompact}. Therefore it is equal to the union of a countable subcollection of previously mentioned balls. 

We have proved that every clopen set is a finite union of base sets. But there are only countably many such unions. \end{proof}

\begin{fact}
The set of all meagre subsets of the topological space $2^\omega$ constitutes an ideal. We call that ideal $\M$.
\end{fact}

\begin{proof} \ 
\begin{enumerate}

    \item The empty set is trivially meagre.
    
    Now we are going to prove that the whole space is not meagre. Take a nowhere dense set $A_0$. Its closure cannot be equal to $2^\omega$ (because its interior is not empty). Therefore there exists $\cl B_{r_0}(x_0)$ disjoint with $A_0$. Now consider another nowhere dense set $A_1$. If its closure contained $\cl B_r(x_0)$, its interior would also have to contain it. Therefore there exists $\cl B_{r_1}(x_1) \subseteq \cl B_r(x_0)$ disjoint with $A_0 \cup A_1$. This is how we inductively construct a countable decreasing sequence of closed sets. Their intersection must contain a point, so we have proved that for any countable sequence of nowhere dense sets we are able to find a point not belonging to their union. Therefore the whole space cannot be meagre. 
    
    \item Let $N, M$ be meagre, $N=\bigcup N_n, M=\bigcup M_n$ and $N_i, M_j$ nowhere dense. Then $N\cup M$ can be written as a countable union of alternatively chosen elements of both sequences.
    
    \item Let $M$ be a meagre set, $N\subseteq M$. Then $M=\bigcup M_n$, where $M_i$ are nowhere dense. But it means that $\int \cl M_n = \emptyset$. Now consider $N_n:=M_n \cap N$. Of course $N_n\subseteq M_n$. But then $\cl N_n\subseteq \cl M_n$, hence $\int \cl N_n\subseteq \int \cl M_n=\emptyset$.
\end{enumerate}
\end{proof}

\begin{df}
The \b{covering coefficient} of an ideal $I$ is the smallest cardinality of a family of sets from an ideal covering the whole ideal. Formally, $\cov(I):=\min\{|{\mathcal A}|:{\mathcal A}\subseteq I \wedge\bigcup{\mathcal A} = X\big\}$.
\end{df}
\begin{fact}
We can use "nowhere dense" in place of "meagre" in the previous definition.
\end{fact}
\begin{proof}
$$2^\omega=\bigcup_{\alpha<\cov(\M)} M_\alpha = \bigcup_{\alpha<\cov(\M)} \bigcup_{n<\omega} N_\alpha^n= \bigcup_{\cone{\alpha,n}\in \cov(\M) \times \omega} N_\alpha^n,  $$
$|\cov(\M) \times \omega|=\cov(\M)$ because of Baire's category theorem, which states that $\cov(\M) > \omega$.
\end{proof}
\begin{thm}
Every space $X$ of weight smaller than $\cov(\M)$ contains a non-trivial convergent sequence. 
\label{covMconv}
\end{thm}
The following proof is a modified version of a proof from \cite{geshke}.
\begin{proof}
We will prove that result for spaces with no isolated points. It will suffice, because either $X$ is scattered, in which case we have already proved that result, or it contains a subspace with no isolated points, in which case we will find a convergent sequence there. Note that the weight of a subspace cannot be bigger than the weight of the original space. 

We are going to construct a continuous surjection $f: X \supseteq Y \to 2^\omega$. To construct it, take the family $\cone{O_s}_{s\in 2^{<\omega}}$ of open subsets of $X$, indexed by finite sequences of zeroes and ones. We also demand that for incomparable $s,t$ the sets $O_s, O_t$ are disjoint and if $s \subseteq t$, then $O_s \subseteq O_t$.  We can easily construct such a family by induction: given $O_s$ we take two different points $x,y \in O_s$ (which we can do by our assumption that no open sets are singletons) and use Hausdorffness to obtain two disjoint open sets separating these to points. After restricting them to $O_s$ we obtain two sets satisfying given conditions. We call them $O_{s\smallfrown 0}, O_{s\smallfrown 1}$.

We now define the subspace $Y:=\bigcap_{n\in \omega} \bigcup _{s\in 2^n}\cl O_s$. Note that any $y\in Y$ is in infinitely many sets $\cl O_s$. But all sequences of equal length are incomparable, so $y$ is in exactly one set $\cl O_t$ if we restrict our attention only to indices of a certain, fixed length.  Of course then $x$ belongs to either $\cl O_{t\smallfrown 0}$ or $\cl  O_{t\smallfrown 1}$.
This means that we can assign a sequence $x$ such that $y \in \cl O_{x|n}$ for all $n\in\omega$. This is the surjection $f$ we were aiming to define. It is continuous, because the preimages of base sets are closed.

We can assume $f$ is irreducible, which means that no proper closed subspace of $Y$ is mapped by $f$ onto the whole $2^\omega$ (because otherwise we can restrict our function to a subspace having that property -- by transfinite induction we obtain a decreasing sequence of such subspaces and take their intersection, which by compactness must be non-empty; it is the subspace we wanted, because for any $x \in 2^\omega$ the set $f^{-1}[\{x\}]$ is closed and so is its intersection with every element of the decreasing sequence -- so the finite subspace cannot be disjoint with this preimage ). For every open $U\subseteq Y$ we define $D_U:=\bigcup\{A\in \Clop(2^\omega): f^{-1}[A] \subseteq U \vee f^{-1}[A] \cap U = \emptyset\}$. 

We are now going to show that every $D_U$ is dense. Let $V$ be an arbitrary non-empty open subset of $2^\omega$. If $f^{-1}[V]$ is disjoint with $U$, then trivially $A\subseteq D_U$. If, however, $f^{-1}[V] \cap U$ is non-empty, it is a non-empty subset of $Y$. 
Note that (from irreducibility) $f$ restricted to the complement of that set is not surjective. Therefore, there exists a non-empty clopen $W\subseteq 2^\omega$ disjoint with $f[Y\setminus (f^{-1}[V]\cap U)$. But then $W\subseteq V$ and $f^{-1}[W] \subseteq U$, therefore $D_U$ is indeed dense.

The set $D_U$ is also open, so its complement $D_U^c$ is nowhere dense. Now we consider the family of nowhere dense sets $D:=\{D_U^c: U \in \mathcal B\}$, where $\mathcal B$ is of cardinality lesser than $\cov(\M)$, so $D$ has to be too. Therefore $D$ cannot cover the whole $Y$, so there must exist $p\in Y\setminus \bigcup D$. 

We will finish the proof by showing that any $y \in f^{-1}[\{p\}]$ has a countable neighbourhood base. It suffices because of Theorem \ref{firstcontconv}. Indeed, we will show that $\{f^{-1}[A]:A \text{ clopen in } 2^\omega \} $ contains a neighbourhood base. Take a point $y \in f^{-1}[\{p\}]$ and its base neighbourhood $U \subseteq Y$. We know that $p\in D_U$, so there is a clopen set $C\ni p$ such that either $f^{-1}[C] \subseteq U$ or $f^{-1}[C] \cap U = \emptyset$. But the second case cannot be true because $y\in f^{-1}[\{p\}] \cap U$. Therefore $f^{-1}[C]\subseteq U$.

\end{proof}

\subsection{Valdivia compacta}
Classes of topological spaces stable under some operations often are subjects of intensive study. For example, \v{C}ech-Stone compactification of $\Sigma$-products is the $\Sigma$-product of \v{C}ech-Stone compactifications. We are going to focus our attention on the $\Sigma$-product of $\R$s. 
\begin{df}
For an arbitrary set $\Gamma$ we define $\Sigma(\Gamma)$ as the set of all the vectors from $\R^\Gamma$ with at most countably many non-zero coordinates.
\end{df}
It turned out that $\Sigma(\Gamma)$ contains all Eberlein spaces and that these spaces are stable under continuous images. This notion has been generalised to Corson spaces, i.e. compact subspaces of $\Sigma(\Gamma)$. They too are stable under continuous images. However, we can generalise Corson spaces, introducing Valdivia compacta. 

The following results can be found in \cite{Kalenda2000} and \cite{Kalenda1999}.
\begin{df}
A $\Sigma$-subset $A$ of a space $K$ is a set for which there exists an injective homeomorphism $h: K \to \R^\Gamma$ such that $h[A]=h[K]\cap \Sigma(\Gamma)$.
\end{df}
\begin{df}
We say that a space $X$ is a \b{Valdivia compactum} if it is compact and contains a dense $\Sigma$-subset.
\end{df}
\begin{df}
A set $K\subseteq X$ is \b{countably compact} if its arbitrary countable open cover has a finite subcover. 
\end{df}
\begin{fact}
A space is countably compact iff every countable decreasing family of closed sets has a non-empty intersection.
\end{fact}
\begin{proof}
If the space is countably compact, the complements of the distinguished sets would be an countable cover without a finite subcover if the intersection was empty.

Conversely, consider a countable open cover $U_0, U_1, U_2, \ldots$. Note that $V_0:=U_0,V_1:=U_0 \cup U_1, V_2:=  U_0 \cup U_1 \cup U_2, \ldots$ is also an open cover. The complements of $V_n$ are decreasing closed sets, therefore there has to be an element outside of all $V_n$. But it would contradict it being a cover. Therefore from some point on, the complements must be empty, so if $V_m^c=\emptyset$, then $X=V_m=U_0\cup\ldots\cup U_m$ is a finite subcover.  
\end{proof}
\begin{df}
A set $F\subseteq X$ is \b{countably closed} if the closure of every countable subset of $F$ is also contained in $F$.
\end{df}
\begin{lm}
\label{sigmactblclosed}
All $\Sigma$-subsets are countably closed.
\end{lm}
\begin{proof}
First, note that $\Sigma(\Gamma)$ is countably closed -- taking closure of a countable set will not add an element with uncountably many non-zero coordinates, because in our set only countably many coordinates are non-zero in at least one element to begin with.

Then notice that countable closedness is transferred back by an injective homeomorphism. 
\end{proof}
\begin{fact}
\label{ctblclosedisctblcompact}
A countably closed subset of a compact space is countably compact.
\end{fact}
\begin{proof}
Take a family $F_n$ of sets closed in a countably closed set $A$. We can assume that the family is decreasing. We need to show it has a non-empty intersection.  Now take $f_n\in F_n$. Consider a family of sets $$ \{ \cl \{f_0,f_1,f_2,\ldots\},\cl \{f_1,f_2,\ldots\},\cl \{f_2,\ldots\},\ldots  \}.$$ They are closed in $A$, so their intersection must be non-empty, so $\bigcap F_n\neq \emptyset$.

\end{proof}
\begin{lm}
\label{urysohn}
Every $\Sigma$-subset $K$ of a compact space has the property that for every $A\subseteq K$ every member of $\cl A$ is a limit point, i.e $\exists A \ni x_n \to x$.
\label{sigmasubsetslimitpoints}
\end{lm}
\begin{proof}
We find the appropriate sequence in $\Sigma(\Gamma)$ and the injective homeomorphism transfers the sequence to $K$.

For $x\in\Sigma(\Gamma)$ at most countably many coordinates of $x$ are non-zero, so we can enumerate them, obtaining a sequence $\cone{\gamma_n(x):n<\omega}$. If it is finite, make it periodic. 

Now take a subset $A\subseteq \Sigma(\Gamma)$ and a vector $x\in\cl A$. We inductively construct a sequence $x_n\in A$ such that $|x_n(\gamma_k(x_l))-x(\gamma_k(x_l))|<\frac1n$ for $0\le l<n, 0\le k \le n$. It is the sequence we wanted. 
\end{proof}

\begin{thm}
An infinite Valdivia compactum contains a non-trivial convergent sequence.
\label{valdiviaconv}
\end{thm}
\begin{proof} 
Let $A$ be a dense $\Sigma$-subset of $K$, an infinite Valdivia compactum. It means that $A$ is infinite. It is also countably compact, by \ref{sigmactblclosed} and \ref{ctblclosedisctblcompact}, so it contains a point $x$ whose all neighbourhoods are infinite (otherwise for every point we take its finite neighbourhood, thus obtaining a countable open cover without a finite subcover). It follows that $x \in \cl (A\setminus\{x\}$ (because $x$ couldn't be in the complement of the closure, because its every neighbourhood would then wander out of the allegedly open complement), so by \ref{urysohn} we obtain the required sequence.

\end{proof}

\section{A space without non-trivial convergent sequences}
\subsection{Construction of \texorpdfstring{$\beta\omega$}{betaomega}}
\begin{df}
A \textbf{filter} over $X$ is a family $F$ of subsets of $X$ such that:
\begin{enumerate}
    \item $\emptyset \not\in F$, $X \in F$;
    \item if $A, B \in F$, then $A \cap B \in F$;
    \item if $A \in F$ and $A\subseteq B$, then $B \in F$.
    \end{enumerate}
\end{df} 
\begin{df}
An \textbf{ultrafilter} over $X$ is a $\subseteq$--maximal filter over $\mathcal{X}$.
\end{df}

\begin{lm}
If $U$ is an ultrafilter over $X$, then for every $A \subseteq X$ either $A \in U$ or ${X} \setminus A \in U$.
\label{eitherinultra}
\end{lm}
\begin{proof}
First notice that at most one of $A, X\setminus A$ can be in a particular filter because a filter is closed under intersections and cannot contain the empty set, but the intersection of $A$ and $X\setminus A$ is empty. 

We will prove that given a set $A \in X$ such that neither it nor its complement is in the filter $F$, it is possible to extend $F$ to a strictly larger filter $G$.
 
Let  $G = \{H\subseteq {X}: (\exists S \in F) S \cap A \subseteq H \}.$  We will show that it is indeed a larger filter, beginning with checking the filter axioms. 

\begin{enumerate}
    \item It is obvious that $X\in G$. If the empty set was in $G$, it would mean that for some $S\in {F}$ the set $S\cap A = \varnothing$, which would then mean that $S \in {X} \setminus A$, so ${X} \setminus A \in {F}$, which is a contradiction. 
    \item If $G_1, G_2 \in {G}$, there exist  $F_1, F_2 \in {F}$ such that $F_1 \cap A \subseteq G_1$ and $F_2 \cap A \subseteq G_2$. Because $F_1 \cap F_2 \in {F}$, then $G_1 \cap G_2 \supseteq (F_1 \cap A) \cap (F_2 \cap A) = (F_1 \cap F_2) \cap A$, hence $G_1 \cap G_2 \in {G}$.
    \item Closure under superset is trivial. 
\end{enumerate}
So far we have shown that $G$ is a filter. What is left is to prove that it is an extension of $F$ containing $A$. 

The set $A$ obviously belongs to $G$, as it suffices to fix $S=X, H=A$ in the definition to obtain $X\cap A \subseteq A$. Now take any $F_1 \in F$ and fix $S=F_1=H$ in the definition. We get $F_1 \cap A \subseteq F_1$. Hence every element of the filter $F$, as well as the set $A$, belongs to the filter $G$ and the proof is complete. 
\end{proof}

Our aim will be to define a topology on the space of all ultrafilters over $\omega$, denoted $\beta\omega$.

\begin{df}
Given any $A\subseteq \omega$ we define a \textbf{cone} of $A$ as the family of all ultrafilters containing $A$, denoted  $\langle A \rangle := \{U \in \beta\omega: A \in U$\}. \end{df}

This definition allows us to concisely describe the simplest ultrafilters -- the principal ultrafilters $U_n := \cone{\{n\}}$.  

\begin{fact} Let us observe three key properties of cones: 
\begin{enumerate}
    \item $\cone{A\cap B} = \cone{A}\cap \cone{B}$,
    \item $\cone{A\cup B} = \cone{A}\cup \cone{B}$,
    \item $\cone{\omega \setminus A} = \beta\omega\setminus \cone{A}$.
\end{enumerate}
\end{fact}
\begin{proof} \ 
\begin{enumerate}
    \item Take $U \in \cone{A \cap B}$. It must then contain $A\cap B$. From the superset property, it must contain $A$ and $B$, so it must be in $\cone{A}$ and $\cone{B}$ and thus in $\cone{A}\cap \cone{B}$. 
    
    Conversely, take $U \in \cone{A}\cap \cone{B}$. It must of course be in $\cone{A}$ and $\cone{B}$, and thus contain $A$ and $B$. From the intersection property $A\cap B \in U$, so $U \in \cone{A \cap B}$. 
    
    \item Take $U \in \cone{A \cup B}$. Therefore $A \cup B \in U$. If $A \in U$ or $B\in U$, then $U\in\cone{A}$ or $U\in\cone{B}$, so $U \in \cone{A} \cup \cone{B}$. Suppose otherwise. Then by Lemma \ref{eitherinultra} we get $\omega \setminus A \in U$ and $\omega \setminus B \in U$. By the intersection property  $(\omega \setminus A) \cap (\omega \setminus B) = \omega \setminus (A \cup B) \in U$, which contradicts our assumption that $A \cup B \in U$, because then $ (\omega \setminus (A \cup B)) \cap (A \cup B) = \emptyset \in U$. 
    
    Conversely, take  $U \in \cone{A} \cup \cone{B}$. Then $U\in\cone{A}$ or $U\in\cone{B}$. Suppose without loss of generality that $U\in\cone{A}$, which means that $A \in U$. Because $A\subseteq A\cup B$, from the superset property we obtain $A \cup B \in U$, hence $U \in \cone{A \cup B}$.
    
    \item Lemma \ref{eitherinultra} lets us conclude that the ultrafilters not contained in $\cone{A}$ (and thus not containing $A$) are strictly these containing $\omega \setminus A$ or, equivalently, strictly these contained in $\cone{\omega\setminus A}$.
\end{enumerate}
\end{proof}
Closure of cones under finite intersections allows us to generate a topology using them as a base. Their closure under complements makes that topology zero-dimensional (which means that it has a clopen base).
\begin{fact}
$\beta\omega$ is Hausdorff.
\end{fact}
\proof We need to prove that every two elements can be separated by open sets. Let us take two different ultrafilters $U, V$. They must differ by an element, let's take such an $A$. Without a loss of generality we have $A\in U$ and $ A \not\in V $. But then it must be that $\omega \setminus A \in V$. It means that $U\in \cone{A}, V \in \cone{\omega \setminus A}$. Therefore by the intersection property $\cone{A} \cap \cone{\omega \setminus A} = \cone{A \cap (\omega \setminus A)} = \cone{\emptyset}=\emptyset$. \qed 
\begin{fact}
$\beta\omega$ is compact.
\end{fact}
\proof Assume otherwise. We must then have an open cover without finite subcover. It can be assumed that this cover consists solely of cones $\{\cone{A_\alpha}:\alpha \in I\}$, because every open set is a union of cones; taking them instead of the original open sets only makes it harder to find a finite subcover. 

For all finite subsets of the index set $I$ we have $\cone{A_{\alpha_0}\cup \ldots \cup A_{\alpha_n} } = \cone{A_{\alpha_0}}\cup \ldots \cup \cone{A_{\alpha_n} } \neq \beta\omega$. So $A_{\alpha_0}\cup \ldots \cup A_{\alpha_n} \neq X$, because $X$ is in all ultrafilters. It follows that  $A_{\alpha_0}^c\cap \ldots \cap A_{\alpha_n}^c \neq \emptyset$, empowering us to generate an ultrafilter with the family $\{A_\alpha^c:\alpha \in I\}$. But then, that ultrafilter is bound to be in the cone of some $A_\beta$, so $A_\beta $ and $A_\beta^c$ are in the same ultrafilter, which is a contradiction. 
\qed
\subsection{Properties of \texorpdfstring{$\beta\omega$}{betaomega}}
We can think that $\omega \subseteq \beta\omega$, equating $n$ to $U_n$. It is worth noting that $\beta\omega$ is, in some sense, the best compactification of $\omega$ (formally, it is the \v{C}ech-Stone compactification of $\omega$). It might seem to be somewhat large for this role, but the following fact will show otherwise. 

\begin{fact}
The set $\omega$ is dense in $\beta\omega$.
\end{fact}
\proof Pick any open set. It must contain a cone $\cone{A}$ because of the definition of the base. Pick any $n\in A$. Then we have $U_n \in \cone{A}$, because $A \in U_n$. \qed
\begin{lm} \label{discreteinhausdorff}
Every infinite Hausdorff space contains an infinite discrete subspace.
\end{lm}
\proof  Let $X$ be an infinite Hausdorff space and $Y$ the set of isolated points. If it's infinite, we are done. If it is not, then $Z:= X \setminus Y \neq \emptyset. $ We will be defining $U_n$ and $p_n$ recursively. Let $p_1 \in Z, U_1=X$. Then we take any unisolated $p_n \in U_{n-1}$ and separate it from $p_{n-1}$ by disjoint open sets $U_n \ni p_n, V_{n-1}\ni p_{n-1}$. If it wasn't possible to find such an unisolated point, that open set would have consisted of $p_k$ and isolated points, but then $p_k$ would have to be isolated itself. 

Then $p_n \not \in \cl \{p_j:j\neq n \}$ because $p_n \in V_n$ and we have $V_n\cap \{p_j:j<n\}\subseteq U_n\cap \{p_j:j<n\}=\emptyset$ while $V_n\cap \{p_j:j>n\}\subseteq V_n\cap U_{n+1}=\emptyset$. Therefore $\{p_n:n\in \omega \}$ is an infinite discrete subspace of $X$.

\begin{thm}
There are no convergent nontrivial sequences in $\beta\omega$.
\end{thm}
\proof Let $V_n$ be a sequence of distinct ultrafilters. Since Hausdorffness is hereditary, $\{V_n:n\in\omega\} $ is infinite and Hausdorff. By the previous lemma we have a subsequence $U_n$, which is a discrete subspace of $\beta\omega$. It means that every ultrafilter has its own cone. Therefore there exists a sequence of subsets of $\omega$ such that $A_n$ is in and only in $V_n$. Define $B_n=A_n \setminus (A_1\cup\ldots\cup A_{n-1})$. That sequence is pairwise disjoint and $B_n$ is in and only in $V_n$ (because $B_n=A_n\cap A_1^c \cap \ldots \cap A_{n-1}^c$). 

If we define $B = \bigcup B_{2n}$ then $B$ is in and only in the even terms of the sequence. But then $\cone{B}$ is an open set containing every other term of the sequence. Therefore our subsequence is divergent, so the original sequence must be too. \qed

\printbibliography
\end{document}